
\input amstex
\documentstyle{amsppt}
\hcorrection{-0.2 cm}
\vcorrection{1.0 cm}

\topmatter

\title Densities of ultraproducts of Boolean algebras
\ifnum\pageno=1
 \footnotemark\ 
\fi
\endtitle   
\author Sabine Koppelberg and Saharon Shelah \endauthor 
\affil  \endaffil
\address \endaddress 
\email sabina\@math.fu-berlin.de, shelah\@math.huji.il  \endemail 
\dedicatory \enddedicatory
\date  \enddate 
\thanks Partially supported by DFG grant Ko 490/7-1 and by the 
Edmund Landau
Center (Jerusalem) for research in Mathematical Analysis, supported by the 
Minerva
Foundation (Germany). In addition, the first author gratefully acknowledges 
the
hospitality of the Department of Mathematics of the Hebrew University of 
Jerusalem
\endthanks

\keywords Boolean algebra, ultraproduct, density, $\pi$-weight 
\endkeywords
\subjclass 03C20, 03E10, 06E05  \endsubjclass
\abstract We answer three problems by J. D. Monk on cardinal 
invariants of Boolean algebras. Two of these are whether taking the 
algebraic density $\pi A$ resp. the topological density $\text{d} A$ of a 
Boolean algebra $A$ commutes with formation of ultraproducts; the 
third one compares the number of endomorphisms and of ideals of a 
Boolean algebra.
\endabstract

\footnotetext{Publication no. 415 of the second author} 

\endtopmatter


\document


In set theoretic topology, considerable effort has been put into the study 
of cardinal invariants of topological spaces, see e.g. \cite{Ju1} and 
\cite{Ho}, \cite{Ju2}. In Monk's book \cite{Mo}, similarly a systematic 
study of cardinal invariants of Boolean algebras is undertaken; in 
particular, the behaviour of these invariants with respect to algebraic 
constructions like taking subalgebras, quotients etc. is investigated. One of 
these is the ultraproduct construction, well known from model theory; 
cf. \cite{ChK}. Many questions on ultraproducts are highly dependent 
on set theory; among the more recent results are those in Shelah' s pcf 
theory dealing with the possible cofinalities $ \text{cf }(\prod_{\alpha 
< \kappa} \lambda_\alpha / D)$ where the $\lambda_\alpha$ are 
regular cardinals, hence well-ordered in a natural way, and the 
ultraproduct has the resulting linear order. 

Monk's book contains a list of 66 problems, three of which are 
answered (consistently) in this paper. 

\medskip

\example{Problem 9} Does there exist a system $(A_i)_{i \in
I}$ of infinite Boolean algebras and an ultrafilter $F$ on $I$ such that 
$\text {d} (\prod_{i \in I} A_i / F) < \vert \prod_{i \in I} \text{d}(A_i) / F 
\vert$? \endexample

\example{Problem 12} Is it true that always $\pi (\prod_{i \in I}
A_i / F) = \vert \prod_{i \in I} \pi (A_i) / F \vert$? \endexample 

\example{Problem 60} Is there a Boolean algebra $A$ such that 
$\vert \text{End } A \vert < \vert \text{Id } A \vert$? \endexample

Here $\pi A$ and $\text{d} A$ are the "algebraic" and the "topological" 
density of $A$, defined by
$$ \gather 
\text{d} A = \text{min } \{ \vert Y \vert: Y \text{ a dense subset of the
 Stone space of } A \} \\
\pi A = \text{min } \{ \vert X \vert: X \text{ a dense subset of } A\} 
\endgather $$
(for more definitions and matters on cardinal functions, see \cite{Mo}). 
Note that we are dealing only with infinite algebras and that, trivially, 
$\omega \leq \text{d} A \leq \pi A$, $\text{d}(\prod_{i \in I} A_i / F) \leq
\vert \prod_{i \in I} \text{d} (A_i) / F \vert$ and $\pi (\prod_{i \in I}
A_i / F) \leq \vert \prod_{i \in I} \pi (A_i) / F \vert$.

In Problem 60, $\text{End } A$ is the set of all endomorphisms, 
$\text{Id } A$ the set of all ideals of $A$.

\medskip

In section 1, we give a positive answer to Problem 12 under SCH. Here 
SCH is the Singular Cardinal Hypothesis: if $2^{\text{cf } \lambda} < 
\lambda$ (so $\lambda$ is singular), then $\lambda^{\text{cf } 
\lambda} = \lambda^+$. However, $\neg$ SCH gives a negative answer to 
both problems 9 and 12:

\proclaim{Theorem A} Assume we have cardinals $\kappa$, $\mu$, and 
$(\lambda_\alpha)_{\alpha < \kappa}$ and an ultrafilter $D$ on 
$\kappa$ such that: $\kappa < \mu = \text{cf } \mu $, $\mu^{< \mu} 
< \lambda_\alpha = \text{cf } \lambda_\alpha$, and the cofinality of 
the ultraproduct $\prod_{\alpha < \kappa} \lambda_\alpha / D$ is less 
than its 
cardinality. Then there is a forcing notion $\Bbb R $ such that 

(a) $\Bbb R$ is $\mu$-complete and satisfies the $(\mu^{< 
\mu})^+$-chain condition; hence forcing with $\Bbb R$ preserves 
all cardinalities and cofinalities outside the interval $[\mu^+, 
\mu^{< \mu})$

(b) for $K \subseteq \Bbb R$ $\Bbb R$-generic over $V$, the 
following holds in $V[K]$:
there are Boolean algebras $(A_\alpha)_{\alpha < \kappa}$ such that 
$\lambda_\alpha = \vert A_\alpha \vert = \pi A_\alpha = 
\text{d} A_\alpha$, but for the ultraproduct $A = \prod_{\alpha < \kappa} 
A_\alpha 
/ D$, 
$$ d(A) 
\leq \pi(A) 
= \text{cf }(\prod_{\alpha < \kappa} \lambda_\alpha / D) 
< \vert \prod_{\alpha < \kappa} \lambda_\alpha / D \vert
= \vert \prod_{\alpha < \kappa} \pi(A_\alpha) / D \vert 
= \vert \prod_{\alpha < \kappa} d(A_\alpha) / D \vert. $$
\endproclaim

Note that SCH is known to be independent from ZFC, modulo some large 
cardinal assumption (see \cite{Ma}). And the assumption of Theorem A 
is a consequence of $\neg$SCH, as follows from pcf theory. A particularly 
easy case is the classical one for $\neg$SCH: assume $\lambda$ is strong 
limit and singular, $\kappa = \text{cf } \lambda $ satisfies 
$2^\kappa < \lambda$, but $\lambda ^\kappa > \lambda^+$; let 
$\mu$ be regular such that $\kappa < \mu < \lambda$. Then there are 
(see \cite{Sh, Ch.II, 1.5}) regular $\lambda_\alpha$ such that 
$\lambda = sup_{\alpha < \kappa} \lambda_\alpha$, $\prod_{\alpha < 
\kappa} \lambda_\alpha / J_\kappa ^{bd}$ has true cofinality 
$\lambda^+$ ($J_\kappa ^{bd}$ the ideal of bounded subsets of 
$\kappa$), hence any uniform ultrafilter $D$ on $\kappa$ gives $ 
\text{cf } (\prod_{\alpha < \kappa} \lambda_\alpha / D) = \lambda^+ 
< \vert \prod_{\alpha < \kappa} \lambda_\alpha / D \vert $.
More generally if $\lambda$ violates SCH, i.e. for 
some $\kappa$, we have $2^\kappa < \lambda$ and $\lambda^\kappa > 
\lambda^+$, 
let $\lambda^\prime $ be 
minimal such that ${\lambda^\prime}^\kappa = \lambda^\kappa$ 
(i.e. ${\lambda^\prime}^\kappa \geq \lambda $); so for every cardinal 
$\rho <
\lambda^\prime$, we have $\rho^\kappa < \lambda^\prime$.
 Now take $\mu = \kappa ^+$ and find, by 
\cite{Sh, Ch.II, 1.5}, an appropriate family 
$(\lambda^\prime_\alpha)_{\alpha < \kappa}$ with limit 
$\lambda^\prime$ and $\text{cf }(\prod_{\alpha < \kappa} 
\lambda^\prime_\alpha / J_{\kappa}^{bd}) = 
{\lambda^\prime}^+$. Moreover we can replace 
${\lambda^\prime}^+$ by any regular cardinal in the interval 
$[{\lambda^\prime}^+, {\lambda^\prime}^\kappa]$; similarly for the strong 
limit 
case; see \cite{Sh, Ch. VIII, \S1}.

Theorem 1.1 below and Theorem A show that the answer to Problem 12 is 
independent from
ZFC. However, it has recently been shown in \cite{RoSh 534, 2.6, 2.7} that 
Problem 9
has a positive answer even in ZFC.

\bigskip

Problem 60 is solved in section 8 by

\proclaim{Theorem B} Assume $\mu$ is a strong limit cardinal satisfying 
$\text{cf } \mu = \omega$ and $2^\mu = \mu^+$. Then there 
is a Boolean algebra $B$ such that $\vert B \vert = \vert \text{End } B
\vert = \mu^+ $ and $\vert \text{Id } B \vert = 2^ {\mu^+}$. 
\endproclaim

The organization of sections 2 to 7 is as follows. In section 2, we 
introduce a first order theory $T$ for Boolean algebras with some 
extra structure which allows (e.g. in ultraproducts $A = \prod_{\alpha < 
\kappa} A_\alpha / D $ of models of $T$) to easily compute $\pi A$. 
In section 3, we construct canonical models $A(p)$ of $T$ from what 
we call valuation functions $p$. In sections 4 to 6, we consider the 
forcing notion $\Bbb P$ of valuation functions, determine its 
completeness and chain conditions, and compute $\text{d} A$ and $\pi A$ 
for 
the canonical algebra $A = A(P)$ constructed from a generic valuation 
function $P$. In section 7, we prove Theorem A. 

For definitions and results on set theory, see \cite{Je}; for Boolean 
algebras, \cite{Ko}.

\bigskip

\beginsection 1. Problem 12 under SCH \par 
\bigskip

We give here a positive answer to Monk\rq s problem 12 
under SCH. Given an ultraproduct $A = \prod_{i \in \kappa} 
A_i/D $ of infinite Boolean algebras, we let $\lambda_i = \pi 
A_i $, so $\omega \leq \lambda_i$. For simplicity of notation, 
we will denote, in this section, by $\prod_{i \in \kappa} \lambda_i/D $ 
both the ultraproduct of the $\lambda_i$ and its cardinality.

Note first that the answer is easy if $\lambda_i \leq 
2^\kappa $ for $D$-almost all $i \in \kappa$ (i.e. if $\{ i 
\in \kappa : \lambda_i \leq 2^\kappa \}$ is in $D$) and 
$D$ is regular. For in this case, each $A_i$ has an infinite 
set of pairwise disjoint elements, so $A$ has cellularity at 
least $ 2^\kappa $ and, on the other hand, $\prod_{i \in 
\kappa} \lambda_i/D \leq 2^\kappa$, hence $ 2^\kappa 
\leq \text{c} A \leq \pi A \leq \prod_{i \in \kappa} 
\lambda_i/D \leq 2^\kappa$. Thus Theorem 1.1 covers the 
interesting case: $ 2^\kappa < \lambda_i $ for $D$-almost 
all $i$. 

\proclaim{1.1 Theorem} (SCH) Assume $ 2^\kappa < 
\lambda_i = \pi A_i $ for all $ i \in \kappa $ and $D$ is an 
ultrafilter on $ \kappa $; let $A = \prod_{i \in \kappa} 
A_i/D$. Then $\pi A = \prod_{i \in \kappa} \lambda_i/D $. 
    \endproclaim

\demo{Proof} We know that $\pi A \leq \prod_{i \in 
\kappa} \lambda_i/D$. Let
$$ \lambda = D-\text{lim } (\lambda_i : i \in \kappa), $$
i.e. $\lambda$ is the least cardinal $\rho$ such that 
$\lambda_i \leq \rho$ holds for all $D$-almost all $i$. 
Without loss of generality, $\lambda_i \leq \lambda $ holds 
for all $i \in \kappa$. \par

\medskip

{\it Claim1. } If $ \theta < \lambda $, then $\theta^\kappa \leq 
\lambda $. \par

To see this, pick $i$ such that $\theta < \lambda_i$. Now if 
$\theta \leq 2^\kappa$, then $\theta^\kappa = 2^\kappa < 
\lambda_i \leq \lambda$. Otherwise, $\kappa < 2^\kappa < 
\theta < \theta^+ \leq \lambda_i$, $(\theta^+)^\kappa = 
\theta^+$ by SCH, so $\theta^\kappa \leq \theta^+ \leq 
\lambda_i \leq \lambda$.

{\it Claim 2. } $\pi A \geq \lambda$. \par

Otherwise pick a dense subset $Y$ of $A$ of size $\rho$, 
where $\rho < \lambda$, say $Y = \{y_\alpha /D : \alpha < 
\rho \} $ with $y_\alpha = (y_\alpha(i))_{i \in \kappa}$ in 
$\prod_{i \in \kappa} A_i $ and $ y_\alpha(i) \neq 0 $. 
Since $\rho < \lambda $, we may assume without loss of 
generality that $\rho < \lambda_i$ for all $i$. So we can 
pick, for $i \in \kappa $, $a_i \in A_i \setminus \{ 0 \}$ 
satisfying $ y_\alpha(i) \nleq a_i $, for all $ \alpha < \rho $. 
The sequence $a = (a_i)_{i \in \kappa} $ is such that 
$y_\alpha / D \nleq a/D $ for $\alpha < \rho $, a 
contradiction.

\medskip

The theorem now follows immediately from the next three 
claims.

\medskip

{\it Claim 3. } If $\pi A \geq \lambda^+$, then the 
assertion of the theorem holds. \par

For in this case, $\lambda^+ \leq \pi A \leq \prod_{i \in 
\kappa} \lambda_i/D \leq \lambda^\kappa / D \leq 
\lambda^\kappa \leq \lambda^+$, where the last inequality 
follows from SCH and $2^\kappa < \lambda$.

{\it Claim 4. } If $\pi A = \lambda$, then every function $ 
f \in \prod_{i \in \kappa} \lambda_i/D $ is bounded below 
$\lambda$, modulo $D$. \par

For the proof, work as in Claim 2: fix a dense subset $Y$ of 
$A$, $ Y = \{y_\alpha /D : \alpha < \lambda \} $, $y_\alpha = 
(y_\alpha(i))_{i \in \kappa}$, $ y_\alpha(i) \neq 0 $. Given 
$f \in \prod_{i \in \kappa} \lambda_i$, we know that $Y_i = 
\{ y_\alpha(i): \alpha < f(i) \}$ cannot be dense in $A_i$, 
since $\vert Y_i \vert \leq \vert f(i) \vert < \lambda_i = 
\pi A_i $. So pick $a = (a_i)_{i \in \kappa}$ where $a_i \in A_i 
\setminus \{ 0\} $ is such that $y_\alpha(i) \nleq a_i $, for 
all $ \alpha < f(i) $. Since $Y$ is dense in $A$, pick $ 
\alpha < \lambda $ such that $y_\alpha / D \leq a / D $. It 
follows that: $y_\alpha(i) \leq a_i$, for $D$-almost all $i$; 
$ \alpha \nless f(i) $ for these $ i $, so $ f(i) \leq \alpha $ ; 
i.e. $ f(i) \leq \alpha $ for $ D $-almost all $ i $. 
Thus $ f $ is bounded by $ \alpha < \lambda $.

{\it Claim 5. } If $\pi A = \lambda$, then the assertion of 
the theorem holds.  \par

For Claim 4 says that for every $f \in \prod_{i \in \kappa} 
\lambda_i$, $f/D = f^\prime / D$ for some $f^\prime : 
\kappa \to \nu $ and some $\nu < \lambda$. By Claim 1, 
$\prod_{i \in \kappa} \lambda_i/D \leq \sum_{\nu < 
\kappa} {\vert \nu \vert}^\kappa \leq \lambda$. It now 
follows from Claim 2 that $ \lambda \leq 
\pi A \leq \prod_{i \in \kappa} \lambda_i/D \leq \lambda 
$.  \qed \enddemo

\bigskip

\beginsection 2. The theory $T$ \par

\bigskip

We sketch here a first order theory $T$. Its relevance for solving 
Problem 12 of \cite{Mo} lies in the fact that the models $\frak A$ of $T$ 
are enlargements $(A, \dots)$ of a Boolean algebra $A$; the extra 
structure of $\frak A$ allows to easily compute $\pi(A)$ --- see 2.1. 
below. Since every ultraproduct $\frak U = (U, \dots)$ of models of $T$ 
is again a model of $T$, we can then similarly compute $\pi(U)$. 

Let $T$ be the first order theory (in an appropriate language) saying 
that, for every model
$\frak A = (A, +, \cdot, -, 0, 1, L, \leq_L, \sim , v, x) $
of $T$, the following hold true. 
 
(a) $ (A, +, \cdot, -, 0, 1) $ is a Boolean algebra.  \par
(b) $L \subseteq A $ is totally ordered by $\leq_L$ and has no greatest 
element. (We do not require any connection between $\leq_L$ and the 
Boolean partial order of $A$, except the one stipulated by (e) below.)  
\par
(c) $v$ is a map from $A$ to $L$; for $l \in L$, 
$ A_l = \{ a \in A: v(a) <_L l \} $
is a subalgebra of $A$. (Hence $(A_l)_{l \in L}$ is an increasing 
sequence of subalgebras of $A$ whose union is $A$.)  \par
(d) $\sim $ is an equivalence relation on $L$ and its equivalence
classes are convex, with respect to $ \leq_L $.  \par
(e) $x$ is a map from $L$ into $A$ (we write $x_i$ for $x(i)$) such 
that $i < l$ implies $x_i \nleq x_l$. Moreover for $l \in L$, the set $ 
\{x_i: i \sim l \}$ is dense for $A_l$ in the sense that for every $a \in 
A_l \setminus \{ 0 \}$ there is some $i \sim l$ satisfying $ 0 < x_i \leq
a$. (Hence $ \{ x_i: i \in L \}$ is a dense subset of $A$.) 
 
\bigskip

\proclaim{2.1 Remark} Let $ \frak A = (A, \dots)$ be a model of $T$, 
$\rho $ the cofinality of the linear order $(L, \leq_L)$ and assume that 
all equivalence classes in $L$ have cardinality at most $\rho$. Then 
$\pi(A) = \rho$. \endproclaim

 \demo{Proof} To see that $\pi(A) \leq \rho$, fix a cofinal subset $M$ 
of $L$ of size $\rho$. The set 
$$ \{x_i : i \sim m, \text{ for some } m \in M \} $$
has size $\rho $ and is dense in $A$, by (e). Assume for contradiction 
that $A$ has a dense subset $X$ of size less than $\rho $. Without loss 
of generality, $X \subseteq \{ x_i: i \in L \} $; pick $l \in L$ such that 
$x_i \in X$ implies $i < l$. $X$ being dense in $A$, there is $x_i \in X$ 
such that $0 < x_i \leq x_l$. So $i < l$ which is impossible by (e). \qed 
\enddemo

\bigskip

In Sections 3 and 4, we will construct "standard" models $ \frak A = (A, 
\dots)$ of $T$ which will roughly look like this, for some regular
cardinal $\lambda$: $\vert A \vert = \lambda$, so without loss of
generality, $\lambda \subseteq A$; we let $L = \lambda $ and $\leq_L$ 
its
natural well-ordering. $A$ will be generated by a sequence $(x_i)_{i \in 
\lambda}$; we then let $A_l$ be the subalgebra of $A$ generated by 
$\{ x_i : i < l \}$ and define $v(a)$ to be the least $i$ such that $a \in 
A_{i+1}$. Finally we will have an infinite cardinal $\mu < \lambda$ and 
define $i \sim l$ iff $i \leq l < i + \mu$ and $ l \leq i < l + \mu $ 
(ordinal addition); thus the equivalence classes will have size $\mu$.
Satisfaction of condition (e) above will be guaranteed by a careful choice of
the generators $x_i$ --- see Proposition 5.1. In particular, $\pi A$ will be 
$\lambda = \vert A \vert.$

\bigskip

\beginsection 3. Valuation functions 

\bigskip

We construct Boolean algebras $A(p)$ from certain functions $p$, the 
so-called valuation functions. Later the Boolean algebras $A(P)$, where 
$P$ will be a generic valuation function, provide the counterexample 
for Problems 9 and 12 in \cite{Mo} looked for.

 We denote the three-element set consisting of the symbols $ \geq , 
\perp , u = \text{ "undefined"}$ by $3$. For any set $w$ with some 
linear order on it (later $w$ will be a subset of some cardinal 
$\lambda$, hence well-ordered), recall that $ [w] ^2 = \{ (i,j) : i < j \text{ 
in } w \}$. 

Given a Boolean algebra $A$ and a family $(x_i)_{i \in w}$ indexed by 
$w$
in $A \setminus \{0 \}$, we can assign to $(x_i)_{i \in w}$ the function 
$p:
[w]^2 \to 3$ defined by 

$$ p(i,j) = \ \geq \text{ if } x_i \geq x_j $$
$$ p(i,j) = \ \perp \text{ if } x_i \perp x_j , \text{ i.e. } x_i \cdot
x_j = 0 $$
$$ p(i,j) = u \text{ otherwise } .$$

Clearly $p$ has then the following properties: 

(1) if $p(i,j) = \ \geq$ and $p(j,k) = \ \geq$ then $p(i,k) = \ \geq$ 
(where $i < j < k$) \par
(2) if $i < j < k$ and $\{ p(i,j), p(i,k) \} = \{ \perp, \geq \}$, then 
$p(j,k) = \ \perp$; similarly if $i < j < k$ and $p(i,j) = \ \perp$, 
$p(j,k) = \ \geq$, then $p(i,k) = \ \perp$. \par 

Let us call a function $p$ satisfying (1) and (2) above a {\it
valuation function } and $w$ its domain.

Conversely, given a valuation function $p: [w]^2 \to 3$, we construct a 
Boolean algebra $A = A(p)$ from $p$ as follows. Let $ \text{Fr } w $ 
be the free Boolean algebra on the set $ \{ u_i : i \in w \}$ of 
independent generators and let $N(p)$ be the ideal in $\text{Fr } w $ 
generated by the elementary products $u_j \cdot u_i$ where $p(i,j) = \
\perp$ resp. $u_j \cdot - u_i$ where $p(i,j) = \ \geq$. Let then $A(p)$ (or 
$A$, for short) be the quotient algebra $\text{Fr } w / N(p) $ and let
$ c: \text{Fr } w \to A(p) $
be the canonical homomorphism. Setting $x_i = c(u_i)$, for $i \in w$, 
we find that the $x_i$ generate $A$. By the very choice of the ideal 
$N(p)$, $p(i,j) = \ \geq$ implies that $x_i \geq x_j$ and $p(i,j) = \ \perp$ 
implies that $x_i \perp x_j$. To see that no other relations than those 
imposed by $p$ hold for the $x_i$, note the following general principle on
construction of Boolean algebras via generators with prescribed relations.

\proclaim{3.1 Remark} Let $E$ be a set of finite partial functions from 
$w$ 
to $\{0, 1 \}$ and let, for $e \in E$, $q_e$ be the elementary product 
$\prod_{e(i) = 1} u_i \cdot \prod_{e(i) = 0} -u_i$ 
 in $\text{Fr } w $. Assume $N$ is the ideal of $\text{Fr } w$ generated 
by the $q_e$, $e \in E $. Then for any function $g : w \to \{ 0, 1 \}$, 
there is an ultrafilter of $\text{Fr } w / N$ including $\{x_i : g(i) = 1\} 
\cup \{-x_i : g(i) = 0\}$ (i.e. $ \{x_i : g(i) = 1\} 
\cup \{-x_i : g(i) = 0\} $ has the finite intersection
property) iff no $e \in E$ is extended by $g$. \endproclaim

This gives the following properties of the $x_i$ in $A = A(p)$, where 
$p$ is a valuation function.

\proclaim{3.2 Remark} $x_i$ is not in the ideal generated by $\{ x_j: j 
> i \}$. In particular, $x_i \neq 0$, the $x_i$ are pairwise distinct, and 
$i < j$ implies that $x_i \nleq x_j$.
\endproclaim
To see this, consider the function $g : w \to \{ 0, 1 \}$ such that $g(k) = 
1$ iff $k = i$ or ($k < i$ and $p(k,i) = \ \geq$). By Remark 3.1,
let $s$ be the ultrafilter of $A$ induced by $g$. Thus $x_i \in s$ but,
for $j > i$, $x_j \notin s$, which shows the claim.

\bigskip

\proclaim{3.3 Remark} $x_i$ is not in the subalgebra of $A$ generated 
by $\{ x_j: j < i \}$. \endproclaim
For consider the functions $g$ and $h$ from $w$ to $\{0,1 \}$ where 
$g$ 
is defined as in the proof of 3.2, $h(k) = g(k)$ for $k \neq i$, but $h(i) =
0$. Let $s$ and $t$ be the corresponding ultrafilters of $A$, $\phi$ 
and $\psi$ the homomorphisms from $A$ to the two-element algebra
corresponding to $s$ and $t$. Now $\phi$ and $\psi$ coincide on $x_j$ 
for all $j < i$, but not on $x_i$.

\bigskip

\beginsection 4. The partial order of valuation functions \par

\bigskip

For the next sections, fix infinite cardinals $\lambda$ and 
$\mu$ such that $\mu^{< \mu} = \mu$, $\mu < \lambda$, and 
$\lambda$ is regular. We shall 
choose $\lambda$ and $\mu$ somewhat more carefully in 
Section 7. Let $\Bbb P ( \lambda, \mu)$ (or $\Bbb P$, for short) be 
the notion of forcing
$$ \Bbb P = \{ p: p \text{ is a valuation function and } \text{dom } 
p \subseteq \lambda \text{ has size less than } \mu \}$$
ordered by reverse inclusion.

\proclaim{4.1 Remark} $\Bbb P$ is $\mu$-closed. \endproclaim

We now build up some machinery for constructing elements of $\Bbb P$ 
with
prescribed properties. Given a set $r$ of relations of the form $x_i \geq
x_j$, $x_i \perp x_j$ (where $i, j \in \lambda$; the relations may be 
thought of
as being formulas in some formal language in the variables $x_i$, $i \in
\lambda$), we define when a relation $\rho$ can be derived from $r$
and we write $r \vdash \rho$: if $\rho$ has the form $x_k \geq x_l$, 
$r \vdash \rho$ iff there are $i_1, \dots, i_m \in \lambda$ such that the
relations $x_k \geq x_{i_1}$, $x_{i_1} \geq x_{i_2}$, ..., $x_{i_m} \geq
x_l$ are all in $r$ (in particular, $ r \vdash x_i \geq x_i $); 
if $\rho$ has
the form $x_k \perp x_l$, $r \vdash \rho$ iff there are $\alpha, \beta \in 
\lambda$ 
such that $x_\alpha \perp x_\beta$ is in $r$ and $r \vdash x_\alpha 
\geq x_k$, $r
\vdash x_\beta \geq x_l$.

Call $r$ {\it consistent } if no relation of the form $x_j \geq x_i$
where $i < j$ and no relation of the form $x_k \perp x_k$ is derivable
from $r$. 
Given $p \in \Bbb P$, define $\text{rel } p$, the relevant part of $p$, 
by
$$ \text{rel } p = \{x_i \geq x_j: p(i,j) = \ \geq \} \cup \{x_i \perp x_j: 
p(i,j) =\ \perp \} .$$ 

\proclaim{4.2 Proposition} If $\vert r \vert < \mu$, then $r$ is consistent
iff $ r \subseteq \text{rel } p$ for some $p \in \Bbb P $.
 
\endproclaim

\demo{Proof} Assume first that $p \in \Bbb P$ and $r \subseteq 
\text{rel } p$ where $\text{dom } p = w \subseteq \lambda $. Then in the
Boolean algebra $A(p)$ constructed in Section 3, all
relations in $r$ and hence all relations derivable from $r$ are satisfied
by the canonical generators $\{x_i: i \in w \}$; moreover, these
generators are non-zero. Hence no relation $x_k \perp x_k$ and no 
relation
of the form $x_j \geq x_i$, $i < j$, can be derived from $r$.

 Conversely, if $r$ is consistent, let $w$ be any subset of $\lambda$ 
such that $\vert w \vert < \mu$ and $\{ i \in \lambda: x_i \text{
occurs in } r \} \subseteq w$. Define $p: {[w]}^2 \to 3$ by

\medskip 

$ p(i,j) = \ \geq \text{ iff } r \vdash x_i \geq x_j$

$ p(i,j) = \ \perp \text{ iff } r \vdash x_i \perp x_j$

$ p(i,j) = u \text{ otherwise}$.

\medskip 

$p$ is a well-defined function (i.e. $r$ does not derive both 
$x_i \geq x_j$ and $x_i \perp x_j$, for $i < j \in w$) since otherwise,
$ r \vdash x_j \perp x_j $, contradicting the consistency of $r$. By the 
above
definition of derivability from $r$, $p$ is a valuation function.
 \qed \enddemo

For further reference, call $p \in \Bbb P$ defined from a consistent set $r$ 
and $w \subseteq \lambda$ as in the proof above the {\it canonical 
extension }
 of $r$ over $w$.

We give one trivial and one not-so-trivial application of this machinery. If 
$G \subseteq \Bbb P$ is $\Bbb P$-generic over our universe $V$ of set 
theory,
then clearly $P_G = \bigcup G$ is a valuation function with $\text{dom } 
P_G =
\bigcup_{p \in G} \text{dom } p$.

\proclaim{4.3 Remark} If $G$ is generic, then $\text{dom } P_G = 
\lambda$.  \endproclaim
To see this, we have to make sure that, for $i \in \lambda$, the set $D_i = 
\{p \in \Bbb P: i \in \text{dom } p \}$ is dense in $\Bbb P$. But given 
$q \in \Bbb P$, let $w \subseteq \lambda$ be such that $\vert w \vert < 
\mu$ 
and $\text{dom } q \cup \{ i \} \subseteq w$. Now by 4.2, $ \text{rel } 
q $ is consistent; let $p$ be the canonical extension of $ \text{rel } 
q $ over $w$. Then $p \in D_i$ and $q \subseteq p$.

\proclaim{4.4 Proposition} If $p, q \in \Bbb P$ coincide on $\text{dom } p
\cap \text{dom } q$, then they are compatible in $\Bbb P$. 
\endproclaim 
\demo{Proof} This follows from a number of claims. We write $p \vdash 
\dots$ 
instead of $\text{rel }p \vdash \dots$ and we say that a relation, e.g. 
$x_i \geq x_j$, is in $p$ if $p(i,j) = \ \geq$ etc.

{\it Claim 1. \/} If $p \vdash x_i \geq x_j$ where $i < j$, then $i, j
\in \text{dom } p$ and the relation $x_i \geq x_j$ is in $p$. Similarly for 
$q$ and for relations of the form $x_i \perp x_j$. --- The claim holds
because $\text{rel } p$, for $p \in \Bbb P$, is closed under
derivations. 

By 4.2 we are left with showing that the set 

$$ r = \text{rel }p \cup \text{rel }q $$

is consistent.  
 
{\it Claim 2. \/} If $r \vdash x_i \geq x_j$, then $p \vdash x_i \geq x_j$ 
 or $q \vdash x_i \geq x_j$ or, for some 
$\alpha$, ($p \vdash x_i \geq x_\alpha$ and $q \vdash x_\alpha \geq 
x_j$) or,
for some $\alpha$, ($q \vdash x_i \geq x_\alpha$ and $p \vdash x_\alpha 
\geq
x_j$).  

{\it Claim 3. \/} If $r \vdash x_i \perp x_j$, then $p \vdash x_i \perp x_j$ 
 or $q \vdash x_i \perp x_j$ or, for some 
$\alpha$, ($p \vdash x_i \perp x_\alpha$ and $q \vdash x_\alpha \geq 
x_j$) or,
for some $\alpha$, ($q \vdash x_i \perp x_\alpha$ and $p \vdash 
x_\alpha
\geq x_j$) (or similarly with $i$ interchanged with $j$).  

{\it Claim 4. \/} If $r \vdash x_i \geq x_j$ and $i, j \in \text{dom } p$,
then $p \vdash x_i \geq x_j$. Similarly for $q$ and for relations of the
form $x_i \perp x_j$. 

The proofs are easy but require consideration of a number of cases. We give 
two
typical examples. In Claim 3, assume e.g. that $p \vdash x_\gamma \perp
x_\delta$, $q \vdash x_\gamma \geq x_i$ and $q \vdash x_\delta \geq 
x_j$.
Then $\gamma$ and $\delta$ are in $ \text{dom } p \cap \text{dom } q 
$, 
$x_\gamma \perp x_\delta$ is (by Claim 1) in $p$, hence in $q$, because 
$p$ and $q$ coincide on $\text{dom } p \cap \text{dom } q$, and $q 
\vdash
x_i \perp x_j$.

Similarly in Claim 4, assume e.g. that $p \vdash x_i \geq
x_\alpha$ and $q \vdash x_\alpha \geq x_j$ where $i, j \in \text{dom } 
p$.
Since $\alpha$ is in $\text{dom } p \cap \text{dom } q$, it follows that 
$x_\alpha \geq x_j$ is in $p$, hence $p \vdash x_i \geq x_j$. 

{\it Claim 5.\/} $r$ is consistent. --- For otherwise by Claim 3, we may
assume that, e.g., for some $\alpha$, $p \vdash x_k \perp x_\alpha$ and 
$q \vdash x_\alpha \geq x_k$. Then $k$ and $\alpha$ are in $\text{dom } 
p
\cap \text{dom } q$, $x_\alpha \geq x_k$ is in $q$ and $x_k \perp x_k$ 
is
in $p$, a contradiction.  \qed 
\enddemo

\bigskip

\proclaim{4.5 Proposition} $\Bbb P$ satisfies the $\mu^+$-chain 
condition. \endproclaim
\demo{Proof} If $X$ is a subset of $\Bbb P$ of size $\mu^+$, then by 
$\mu^{< \mu} = \mu$ and the $\Delta$-lemma there are $p$ and $q$ 
in $X$ coinciding on $\text{dom } p \cap \text{dom }q$. So we are 
finished by Proposition 4.4.  \qed \enddemo

\bigskip

\beginsection 5. Computing $\pi (A(P))$ \par \bigskip

In this and the following section, let $G$ be a $\Bbb P$-generic filter over 
$V$ and $P$ the resulting generic valuation function (see 4.3). Write $A$ 
for $A(P)$. We prove condition (e) of section 2 for $A$, thus being able to 
compute $\pi (A)$ in $V[G]$.

\proclaim{5.1 Proposition} The following holds in $V[G]$. Let $\alpha < 
\lambda$
be an ordinal, $a \subseteq \alpha$ finite, $e: a \to \{0,1 \}$ and 
$$ y = \prod_{e(i) = 1} x_i \cdot \prod_{e(i) = 0} -x_i > 0 \qquad (\text{in } 
A ).$$
Then there is $i^* \in [\alpha, \alpha + \mu)$ (ordinal addition) 
such that $x_{i^*} \leq y$. - In particular, the set $\{ x_{i^*} :i^* \in 
[\alpha, \alpha + \mu) \}$ is dense for the subalgebra of $A$ generated by 
$\{x_i : i < \alpha \}$.   \endproclaim

\demo{Proof} We do not distinguish notationally between elements of 
$V[G]$ 
and their $\Bbb P$-names; in particular since $a$ and $e$, being finite, are 
in the ground model. Pick $p \in G$ such that
$$ p \Vdash y = \prod_{e(i) = 1} x_i \cdot \prod_{e(i) = 0} -x_i > 0; $$
it suffices to prove that
$$ D = \{t \in \Bbb P :t \leq p, \text{ and } t \Vdash x_{i^*}\leq y \text{ for 
some }
 i^* \in [\alpha, \alpha + \mu) \} $$
is dense below $p$. To this end, let $q \leq p$ be arbitrary. By 4.3, we can 
fix 
$r \leq q$ such that $a \subseteq \text{dom } r$. Then fix $i^* \in
[\alpha, \alpha + \mu) \setminus \text{dom } r$;
this is possible by $\vert \text{dom } r\vert < \mu$. We define a function 
$s$
with domain $a \cup \{i^* \}$ by putting
$$
\gather s \restriction {[a]}^2 = r \restriction {[a]}^2 \\
s(i, i^*) = \ \geq \text{ if } i \in a \text{ and } e(i) = 1 \\
s(i, i^*) = \ \perp \text{ if } i \in a \text{ and } e(i) = 0.
\endgather 
$$ 

{\it Claim} $s \in \Bbb P$, i.e. $s$ is a valuation function. 

Let us check just one case. Note that, for $u \in \Bbb P$, 
$u(i,j) = \ \geq $ implies that $u \Vdash x_i \geq x_j$ and similarly 
for $\perp$ instead of $\geq$ since for any generic $H \subseteq \Bbb P$ 
containing $u$, $u \subseteq P_H$ and thus $x_i \geq x_j$ will hold in 
$A(P_H)$. Assume e.g. $i < j$ in $a$, $s(i, j) = \ \geq $ and $s(j, i^*) = \ 
\geq$;
we have to show that $s(i, i^*) = \ \geq$. 
The assumptions say that $r(i, j) = \ \geq$ (since $i, j \in a$) and $e(j) = 
1$;
we have to show that $e(i) = 1$. But if $e(i) = 0$, then: 
$p \Vdash 0 \neq -x_i \cdot x_j$ (because $p \Vdash 0 < y \leq -x_i \cdot 
x_j$),
$r \Vdash 0 \neq -x_i \cdot x_j$ (since $r \leq p$), $r \Vdash x_i \geq 
x_j$ 
(by the above assumption), $ r \Vdash -x_i \cdot x_j = 0 $, a contradiction.
Now $r$ and $s$ coincide on $a = \text{dom } r \cap \text{dom } s $, so by 
4.4, pick $t \in \Bbb P$ extending both $r$ and $s$. Then 
$t \leq q$ and $s \Vdash x_{i^*} \leq y$, by the very definition of 
$s$ above, so $t \in D$.
\qed \enddemo

\proclaim{5.2 Corollary} $\pi(A) = \lambda$ (in $V[G]$). 
\endproclaim

\demo{Proof} This follows from Remark 2.1 and the sketch of the model
$\frak A = (A, \dots) \models T$ following it, plus 5.1. 
Let us remark that 6.1 gives another proof, since $\text{d} A = \lambda $, 
 $\text{d} A \leq \pi A$ holds for all Boolean algebras and 
 $\pi A \leq \vert A \vert = \lambda$. 
\qed \enddemo

\example{5.3 Example} Our construction of $A = A(P)$ and 5.1 above give a 
counterexample to the 
assertion in 4.1 of \cite{Mo}, in $V[G]$. For this, let $A_\alpha$ be the 
subalgebra 
of $A$ generated by $ \{ x_i: i < \alpha \}$; so if 
$\alpha \in I =\{\alpha < \lambda: \text{cf } \alpha = \mu \} $, then by 
Remark 2.1 and 5.1 above, we have $\pi A_\alpha = \mu$. Moreover 
$A = \bigcup_{\alpha \in I} A_\alpha $ and $\pi A = \lambda$ where 
$\lambda$ can be larger than $\mu^+$. --- In fact, the argument given in 
\cite{Mo,
4.1} depends on the assumption that the chain $(A_\alpha)_{\alpha \in I}$ 
is 
continuous which is not the case here.
\endexample 

\bigskip

\beginsection 6. Computing $\text{d} (A(P))$ \par

\bigskip

Our single theorem here is the following.

\proclaim{6.1 Theorem} 	In $V[G]$, $A = A(P)$ satisfies $\text{d}(A)
= \lambda$.  \endproclaim

\demo{Proof} Otherwise, the cardinal $\sigma = {\text{d} (A)}^{V[G]}$ is 
less
than $\lambda$. There are a $\Bbb P$-name $u$ and a condition $p \in 
\Bbb
P$ (in fact, $p \in G$) such that
$$ p \Vdash u \text{ is a sequence } (u_\nu)_{\nu < \sigma} \text{, each
} u_\nu \text{ is an ultrafilter of } A \text{, and } A \setminus \{ 0 \} =
\bigcup_{\nu < \sigma} u_\nu . $$
For $\alpha < \lambda$, fix $p_\alpha \in \Bbb P$ and $\nu_\alpha < 
\sigma$ 
such that $p_\alpha \leq p$ and 
$$ p_\alpha \Vdash x_\alpha \in u_{\nu_\alpha} $$
($x_\alpha$ the (name of the) $\alpha$\rq th generator of $A$). In the 
next
steps, we construct stationary subsets $S_1 \supseteq S_2 \supseteq S_3
\supseteq S_4$ of $\lambda$.

{\it Step 1. } $S_1 = \{\alpha \in \lambda : \text{cf } \alpha = \mu \}$ is
stationary in $\lambda $ because $\mu < \lambda$ and $\lambda$ is 
regular.

{\it Step 2. } Since $\sigma < \lambda = \text{cf } \lambda $, there are 
$\nu^* < \sigma$ and a stationary $ S_2 \subseteq S_1 $ such that 
$\nu_\alpha = \nu^*$, for all $\alpha \in S_2$.

{\it Step 3. } Write $w_\alpha = \text{dom } p_\alpha$ , for $\alpha \in
\lambda$. We find $\alpha^* \in \lambda$ and a stationary $ S_3 
\subseteq
S_2$ such that for all $\alpha \in S_3$, $\alpha^* < \alpha$ and 
$w_\alpha \cap \alpha \subseteq \alpha^*$ hold. To this end, note that 
$\text{cf }\alpha = \mu$ for $\alpha \in S_2$ and $ \vert w_\alpha \cap
\alpha \vert < \mu$; so pick $j_\alpha < \alpha$ satisfying $w_\alpha \cap 
\alpha \subseteq j_\alpha$. Apply Fodor\rq s theorem to obtain $S_3$.

{\it Step 4. } We find a stationary set $ S_4 \subseteq S_3$ such that 
$\alpha < \beta$ in $S_4$ implies $w_\alpha \subseteq \beta$. To do 
this,
define by induction $f: \lambda \to \lambda$ strictly increasing and 
continuous
such that, for all $\alpha$, $\bigcup_{\nu < \alpha} w_\nu \subseteq
f(\alpha)$ and let $S_4 = S_3 \cap C$ where $C = \{ \alpha : f(\alpha)
= \alpha \}$ is closed unbounded. Then $S_4$ is stationary and, for 
$\alpha 
< \beta$ in $S_4$, we have $w_\alpha \subseteq f(\beta) = \beta$.

Now $\mu^+ \leq \lambda$ and $\Bbb P$ satisfies the $\mu^+$-chain 
condition. So we can find $\alpha < \beta$ in $S_4$ such that $p_\alpha$ 
and $p_\beta$ are compatible in $\Bbb P$. Let $r$ be the following set 
of
relations:
$$ r = \text{rel }(p_\alpha) \cup \text{rel }(p_\beta) \cup \{ x_\beta 
\perp x_\alpha \}$$
(see the machinery in section 4). 

{\it Claim. } $r$ is consistent.

By the claim and 4.2, pick then $q \in \Bbb P$ such that $r \subseteq
\text{rel}(q)$. This $q$ will force the following statements:  \par
$x_\beta \perp x_\alpha $ \par
$ x_\alpha \in u_{\nu_\alpha} = u_{\nu^* } \text{ and } 
   x_\beta \in u_{\nu_\beta} = u_{\nu^* }$ \par
$u_{\nu^*}$ has the finite intersection property (being an ultrafilter), \par
and this contradiction finishes the proof.

{\it Proof of the Claim. } Clearly no relation $x_i \geq x_j$ where $j < i$ 
can have a derivation from $r$, since such a derivation would not use the
relation $ x_\beta \perp x_\alpha $; hence $x_i \geq x_j$ would be 
derivable 
from $ \text{rel }(p_\alpha) \cup \text{rel }(p_\beta) $, contradicting the
compatibility of $p_\alpha$ and $p_\beta$. \par
Now assume $r \vdash x_k \perp
x_k$, for some $k \in \lambda$. A derivation witnessing this starts, without
loss of generality, with the relation $ x_\beta \perp x_\alpha$. So in 
$p_\alpha \cup p_\beta$ there are relations
$$ x_{i_0} \geq x_{i_1}, \dots, x_{i_{r-1}} \geq x_{i_r} \text{ where } 
     i_0 = \alpha, i_r = k $$
$$ x_{j_0} \geq x_{j_1}, \dots, x_{j_{s-1}} \geq x_{j_s} \text{ where } 
     j_0 = \beta, j_s = k .$$
Note that $\alpha = i_0 < i_1 < \dots < i_r = k$ (since if $x_j \geq x_i$ is 
in $p_\alpha \cup p_\beta$, then $j < i$);
similarly, $\beta = j_0 < j_1 < \dots < j_s = k$ . \par
We prove by induction on $t \in \{0, \dots, r \}$ that $i_t \notin w_\beta
= \text{dom } p_\beta$; for $t = r$ this gives a contradiction because then 
$k = i_r \notin w_\beta$, so $k \in w_\alpha$ and $k \geq \beta$, but 
$w_\alpha \subseteq \beta$. First, $i_0 \notin w_\beta$: otherwise, by Step 
3, 
$i_0 = \alpha \in w_\beta \cap \beta \subseteq \alpha^*$, contradicting 
$\alpha^* < \alpha$ for $\alpha \in S_3$. If $i_t \notin w_\beta$ but 
$i_{t+1} \in w_\beta$, then the relation $x_{i_t} \geq x_{i_{t+1}}$ 
must be in $p_\alpha$. But then $i_{t+1} \in w_\alpha \subseteq \beta$ 
and
again $i_{t+1} \in w_\beta \cap \beta \subseteq \alpha^* < \alpha$, a
contradiction. \qed  \enddemo

\bigskip

\beginsection 7. Proof of Theorem A \par \bigskip

\demo{7.1 Proof of Theorem A} Fix $\kappa$, $\mu$, $\lambda_{\alpha}$ 
and $D$ 
as given in the theorem; $\Bbb R$ will be the iteration of two 
forcing notions. In the first step, collapse $\mu^{ < \mu}$ to 
$\mu$ with $\Bbb Q = \text{Fn } (\mu, \mu^{ < \mu}, < \mu ) $ 
in Kunen\rq s notation (\cite{Ku}). This forcing is 
$\mu$-closed and satisfies the ${(\mu^{ < \mu})}^+$-chain 
condition; in the resulting generic model $V[H]$, $\mu^{ < 
\mu} = \mu$ holds. The notions of ultrafilters on $\kappa$, 
the cartesian product $\prod_{\alpha < \kappa} 
\lambda_\alpha$ etc. are absolute for this forcing by 
$\mu$-closedness of $\Bbb Q$ and $\kappa < \mu$; thus all 
assumptions of the theorem continue to hold in $V[H]$.

Working now in $V[H]$, let, for $\alpha \in \kappa$, $ \Bbb 
P_\alpha$ be the forcing notion $ \Bbb P (\lambda_\alpha, \mu)$ defined
in section 4; let $\Bbb P$ be the full cartesian product
$ \Bbb P = \prod_{\alpha < \kappa} \Bbb P_\alpha $
with the coordinate-wise partial order. For $G \subseteq 
\Bbb P$ generic over $V$, $G_\alpha = 
{\text{pr}_\alpha}^{-1} [G] $ is $\Bbb P_{\alpha}$-
generic over $V[H]$ ($\text{pr}_{\alpha}$ the 
$\alpha$\rq th projection). $\Bbb P$ is clearly $\mu$-
closed, moreover, as in the proof of 4.5, the $\Delta$-lemma 
implies that $\Bbb P$ satisfies the $\mu^+$-chain condition 
since $\mu^{ < \mu} = \mu$. Thus the assumptions of the 
theorem, as well as $\mu^{ < \mu} = \mu$, continue to hold in 
$V[H][G]$.

In $V[H][G]$, $P_{\alpha} = \bigcup G_{\alpha} : 
[\lambda_{\alpha}] ^2 \to 3 $ is a generic valuation 
function. Let $ A_{\alpha} = A(P_{\alpha}) $ be its associated 
Boolean algebra; by sections 5 and 6, $\pi(A_{\alpha}) = 
\text{d} (A_{\alpha}) = \lambda_{\alpha}$. In the standard 
model $\frak A_\alpha = (A_\alpha, \dots) $ of $T$ (see 
section 2), the predicate $L$ is interpreted by 
$\lambda_{\alpha}$ and the equivalence classes of $\sim_L$ 
have size $\mu$. So in the ultraproduct $\frak A = 
\prod_{\alpha < \kappa} \frak A_{\alpha } / D$, $L$ is 
interpreted by $ \prod_{\alpha < \kappa} \lambda_{\alpha } 
/ D$ and the equivalence classes of $\sim_L$ have size $ 
\leq \vert \mu^\kappa / D \vert = \mu$ (by $\kappa < 
\mu $ and $\mu^{ < \mu} = \mu$). Now Remark 2.1 says that 
$\pi(A) = \text{cf }\prod_{\alpha < \kappa} \lambda_{\alpha 
} / D $ and hence
$ \text{d}(A) \leq \pi(A) 
= \text{cf } (\prod_{\alpha < \kappa} \lambda_{\alpha } / 
D) 
 < \vert \prod_{\alpha < \kappa} \lambda_{\alpha } / D 
\vert 
 = \vert \prod_{\alpha < \kappa} \pi(A_{\alpha } )/ D \vert 
= \vert \prod_{\alpha < \kappa} \text{d}(A_{\alpha } )/ D \vert 
 $.  \qed   \enddemo

We can prove a little more:

\proclaim{7.2 Remark} In $V[H][G]$, let $ A = \prod_{\alpha < \kappa} 
A_\alpha /
D$ be the algebra constructed in 7.1 and let $ \lambda = 
\text{cf } \prod_{\alpha < \kappa} \lambda_\alpha 
 / D $. Then $\text{d}(A) = \lambda$. \endproclaim

\demo{Proof} Our proof will closely follow that of 6.1.

Fix a sequence $(f_\gamma)_{\gamma \in \lambda}$ in 
$\prod_{\alpha < \kappa} \lambda_\alpha $ such that 
$(f_\gamma / D)_{\gamma \in \lambda}$ 
is strictly increasing and cofinal in the ultraproduct $\prod_{\alpha <
\kappa} \lambda_\alpha / D$. By \cite{Sh, Ch.II}, the set 
$$
\split 
S = \{ \gamma \in \lambda : &\text{ cf } \gamma = \mu^+ \text{ , and there
is } g \in \prod_{\alpha <\kappa} \lambda_\alpha \\ 
&\text{ such that } g/D 
\text{ is the least upper bound of } \{f_\delta / D :\delta < \gamma \}\\ 
&\text{ and } \text{cf } g(\alpha) = \mu^+ \text{ for all }\alpha \in 
\kappa \}
\endsplit
$$ 
is stationary; so we may assume that, for $\gamma \in S$, $f_\gamma$ 
satisfies the requirements for $g$ above.

Now note that, in $V[H][G]$, $ \text{d} A \leq \pi A = \lambda $ 
as shown in the proof of 7.1; so assume for contradiction that 
$\text{d}A < \lambda$. Thus , in $V[H][G]$,there are a $\Bbb P$-name $u$, 
$\sigma < \lambda$ and $p \in \Bbb P$ such that
$$ p \Vdash u = (u_\nu)_{\nu < \sigma} \text{ is a sequence of ultrafilters
 of } A \text{ covering } A \setminus \{ 0 \} . $$
For $\gamma \in S$, fix $p_\gamma \geq p$ and $\nu_\gamma \in 
\sigma$ such that 
$$ p_\gamma \Vdash y_\gamma / D \in u_{\nu_\gamma} $$
where $y_\gamma $ is (a $\Bbb P$-name for) 
$(x_{f_\gamma(\alpha)})_{\alpha < \kappa} / D $ and $x_i$ 
is (a $\Bbb P$-name for) the $i$\rq th canonical generator of $A_\alpha$, 
for $i < \lambda_\alpha$. There is a stationary subset $S_1$ of $S$ 
such that $\nu_\gamma $ is some fixed $\nu^*$, for $\gamma \in S_1$ 
(because $\nu_\gamma < \sigma < \lambda$ and $\lambda$ is regular).
As in Step 3 in the proof of 6.1, there exists , for $\gamma \in S_1$, some 
$\beta_\gamma < \gamma$ such that, for $D$-almost all $\alpha$,
$$\text{dom } p_\gamma(\alpha) \cap f_\gamma(\alpha) \subseteq 
f_{\beta_\gamma}(\alpha) . $$
Without loss of generality (i.e. by passing to a stationary subset), 
$\beta_\gamma $ is some fixed $\beta^*$, for all $\gamma \in S_1$. Now
$K_\gamma = \{\alpha \in \kappa : \text{dom } p_\gamma(\alpha) \cap
f_\gamma(\alpha) \subseteq f_{\beta^*}(\alpha) \} \in D $, for $\gamma 
\in S_1$; 
since $ 2^\kappa < \lambda$, we may assume without loss of generality that 
$K_\gamma $ is some fixed $K^* \in D$, for $\gamma \in S_1$. 

As in Step 4 of the proof of 6.1, we may assume that $\gamma < \delta$ in 
$S_1$ 
implies that
$$ K_{\gamma \delta} = \{\alpha \in \kappa : \text{dom } 
p_\gamma(\alpha) \subseteq
f_\delta(\alpha) \} \in D $$ 
because $(f_\delta / D)_{\delta \in \lambda}$ is cofinal in $\prod_{\alpha
< \kappa} \lambda_\alpha / D $.

Now $\Bbb P$ satisfies the $\mu^+$-chain condition and $S_1$ has size 
$\lambda \geq \mu^+$; so fix $\gamma < \delta$ in $S_1$ such that
$p_\gamma$ and $p_\delta$ are compatible in $\Bbb P = 
\prod_{\alpha \in \kappa} \Bbb P_\alpha$, i.e. $p_\gamma(\alpha)$ and 
$p_\delta(\alpha)$ are compatible in $\Bbb P_\alpha$, for all $\alpha \in 
\kappa$.

We conclude as in 6.1: for all $\alpha \in K^* \cap K_{\gamma \delta}$, the 
set
$$ r_\alpha = \text{rel } p_\gamma(\alpha) \cup 
\text{rel } p_\delta(\alpha) \cup \{ x_{f_\delta(\alpha)} \perp 
x_{f_\gamma(\alpha)} \} $$
is consistent; so pick $q_\alpha \in \Bbb P_\alpha$ satisfying 
$r_\alpha \subseteq \text{rel } q_\alpha$. Choose $q \in \Bbb P$ having 
$\alpha$\rq th coordinate $q_\alpha$, for $\alpha \in K^* \cap 
K_{\gamma \delta}$; 
then $q$ forces that: 
$y_\delta / D \perp y_\gamma / D $, 
$y_\gamma / D \in u_{\nu_\gamma} = u_{\nu^*} $ and 
$ y_\delta / D \in u_{\nu_\delta} = u_{\nu^*}$, $u_{\nu^*}$ is an ultrafilter. 
This gives a contradiction. \qed

\enddemo

\bigskip

\beginsection 8. Proof of Theorem B   \par

\bigskip
 
To abbreviate the main body of the proof, we state in advance two easy 
lemmas. The proofs are left to the reader.

\proclaim {8.1 Lemma} Assume $h: C \to D$ is a homomorphism of 
Boolean 
algebras, $\{ c_n: n \in \omega \}$ is a partition of unity in $C$, and also 
$\{ h(c_n): n \in \omega \}$ is a partition of unity in $D$. Then, if $x_n 
\in C$ are such that $\sum^C_{n \in \omega} x_n \cdot c_n$ exists, we 
have $ h(\sum^C_{n \in \omega} x_n \cdot c_n ) = \sum^D_{n \in \omega} 
h(x_n \cdot c_n)$.  \endproclaim

Given a subalgebra $C$ of $D$ and $x \in D$, let $I_C(x) = \{ c \in C: c 
\cdot x = 0 \}$, an ideal of $C$. Call $x, y \in D$ equivalent over $C$ 
(and write $x {\sim}_C y$) if both $I_C(x) = I_C(y)$ and $I_C(-x) = 
I_C(-y)$ hold, i.e. if $x$ and $y$ realize the same quantifier-free type 
over $C$.

\proclaim {8.2 Lemma} If $x, y \in D$ are equivalent over $C$, then there 
is
no $c \in C \setminus \{ 0 \}$ disjoint from $x + -y$. \endproclaim

We break up the proof of Theorem B into eight preparatory steps in which 
certain objects are constructed or notation is fixed, plus four claims. Let $C 
\leq D $ denote that $C$ is a subalgebra of $D$; $\overline A$ is the 
completion of $A$.

\medskip

{\it Step 1.} Take $\mu$ as assumed in the theorem, fix a set $U$ of 
cardinality $\mu$, and let $A = \text{Fr } U$, the free Boolean algebra over 
$U$. Since $\vert \overline A \vert = \mu^\omega \geq \mu^+ = 2^\mu$, 
we have 
$\vert \overline A \vert = \mu^+$. The algebra $B$ promised in the 
theorem 
will be a subalgebra of $ \overline A $, generated by $A$ and
pairwise distinct elements $b_i$ of $ \overline A $, $i < \mu^+$.
So $\vert B \vert = \mu^+$ and we know in advance that 
$\mu^+ \leq \vert \text{End } B \vert$ and $\vert \text{Id }B \vert \leq
2^{\mu^+}$.

\medskip

{\it Step 2.} Fix an enumeration $\{ g_j: j < \mu^+ \}$ of all 
homomorphisms 
from $A$ into $\overline A $. This is possible since $\vert A \vert = 
\mu$ 
and $\vert \overline A \vert = \mu^+ = (\mu^+)^\mu $.

\medskip

{\it Step 3.} Fix a sequence $(\mu_n)_{n \in \omega}$ of cardinals such 
that 
$\mu = sup_{n \in \omega} \mu_n$ and $2^{\mu_n} < \mu_{n+1}$.

\medskip

{\it Step 4.} For each ordinal $i < \mu^+$, fix subsets $S_{in}$ of $i$ 
such that $i = \bigcup_{n \in \omega} S_{in}$, $S_{in} \subseteq S_{i,n+1}$ 
and $\vert S_{in} \vert \leq \mu_n$. This is possible since $\vert i 
\vert \leq \mu$.

\medskip

{\it Step 5.} Fix a sequence $(A_n)_{n \in \omega}$ of subalgebras of $A$ 
such that $A = \bigcup_{n \in \omega} A_n$, $A_n \subseteq A_{n+1}$ 
and $\vert A_n \vert \leq \mu_n$. 

\medskip

{\it Step 6.} Define a tree $T = \bigcup_{n \in \omega} T_n$ with $n$'th
level $T_n = \mu_0 \times \dots \times \mu_{n-1}$ where $t \leq s$ in 
$T$ 
means that $s$ extends $t$; so $\vert T \vert = \mu$. The cartesian 
product $F = \prod_{n \in \omega} \mu_n$ has size $\mu^\omega = 
\mu^+$; fix a one-one enumeration $\{f_i: i < \mu^+ \}$ of $F$.

Split $U 
\subseteq A = \text{Fr } U$ (cf. Step 1) into two disjoint subsets $X$ and 
$Z$ 
such that $\vert X \vert = \vert Z \vert = \mu$; then split both $X$ and 
$Z$ 
into pairwise disjoint subsets $X_t$, $t \in T$, and $Z_t$, $t \in T$, such 
that 
$\vert X_t \vert = \mu$ and $Z_t \neq \emptyset$.

\medskip

{\it Step 7.} Here we define, for $i < \mu^+$, the elements $b_i$ of 
$\overline A$ and then let $B$ be the subalgebra of $\overline A$ 
generated by $A
\cup \{b_i : i < \mu^+ \}$. $b_i$ is constructed out of certain
elements $x_{in}$, $y_{in}$, $z_{in}$, $n \in \omega$, of $U$ by putting 
$$ s_{in} = x_{in} + -y_{in} $$ 
 $$ d_{i,n} = s_{i,n} \cdot \prod_{m < n} -s_{im} $$
$$ b_i = \sum_{n \in \omega} z_{in} \cdot d_{in} .$$ 
To choose the $x_{in}$, $y_{in}$, $z_{in}$, fix $i < \mu^+$ and $n \in 
\omega$; thus
$ t = f_i \restriction n $ 
is an element of the tree $T$. Pick $z_{in} \in Z_t$ (see Step 6) 
arbitrarily. $x_{in}$ and $y_{in}$ are chosen much more carefully: we 
want them to be distinct elements of $X_t$ satisfying 

\medskip

(*) for all $j \in S_{in}$, $g_j(x_{in}) {\sim}_{A_n} g_j(y_{in})$ 

\medskip

(cf. Steps 4, 2, 5, and the definition of ${\sim}_{A_n}$ before 8.2). This
is possible since: \par
$\vert A_n \vert \leq \mu_n$ \par
there are at most $2^{\mu_n}$ equivalence classes in $\overline A$, with 
respect 
to ${\sim}_{A_n}$, since there are at most $2^{\mu_n}$ ideals in $A_n$ 
\par
$\vert S_{in} \vert \leq \mu_n$ \par
the set $\{ (g_j(x)/ {\sim}_{A_n})_{j \in S_{in}} : x \in X_t \}$ has size 
at most $2^{\mu_n}$ \par
$2^{\mu_n} < \mu = \vert X_t \vert $ . 

\medskip

{\it Step 8.} (Remark) For $b \in A$, let us denote by $\text{supp } b$ 
(the
support of $b$) the smallest subset of $U$ generating $b$.
Now for $i < \mu^+$, the supports $ \{\text{supp } s_{in} : n \in 
\omega \}$ are
pairwise disjoint and thus $\sum^{\overline A} s_{in}= 1$. It follows that 
the
pairwise disjoint set $\{d_{in} : n \in \omega \}$ is a partition of unity in 
$\overline A$ and all $d_{in}$ are non-zero. --- Similarly, for any 
homomorphism 
$g: A \to \overline A $, the sets $\{g( d_{in}): n \in \omega \}$ and $\{g(
s_{in}): n \in \omega \}$ have the same upper bounds in $A$ resp. 
$\overline A$.

\medskip

{\it Claim 1. } If $j < i < \mu^+$, then $ \{g_j( d_{in}): n \in \omega \}$ 
is a partition of unity (in $\overline A$). --- Otherwise, assume $a \in A^+$ 
and $a \cdot g_j( s_{in}) = 0$ for all $n$ (cf. Step 8). Pick $n$ so large
that $a \in A_n$ and $j \in S_{in}$. Then $a \cdot g_j(x_{in} + -y_{in}) =
0$, so $a \cdot (g_j(x_{in}) + -g_j(y_{in})) = 0$, contradicting (*) and
Lemma 8.2.

\medskip

{\it Claim 2. } Let $g$ be an endomorphism of $B$, say $ g \restriction A 
=
g_j$ (see Step 2). Then for all $i > j$, $ g(b_i) = \sum^{\overline A}
g_j(z_{in}) \cdot g_j(d_{in}) $ holds. Hence $g$ is uniquely determined by
its action on $ A \cup \{b_i: i \leq j \}$. --- This follows from Claim 1 and
Lemma 8.1.  \medskip

{\it Claim 3. } $ \vert \text{End } B \vert \leq \mu^+$. --- To completely 
describe
some $g \in \text{End } B$, we have only $\mu^+$ choices for $g 
\restriction A$ 
(Step 2) and, for $j < \mu^+ $, at most $(\mu^+)^{\vert j \vert} \leq 2^\mu 
=
\mu^+$ choices for $(g(b_i))_{i \leq j}$, so we are finished by Claim 2.

\medskip

{\it Claim 4. } The generators $\{ b_i: i < \mu^+ \}$ are
ideal-independent; hence $\vert \text{Id } B \vert = 2^{\mu^+} $. ---
We prove that, for $i \in \mu^+$ and $J$ a finite subset of $\mu^+ 
 \setminus \{ i \}$, $b_i \nleq \sum_{j \in J} b_j$. (It follows that the 
ideals $I_K$ generated by $\{b_i: i \in K \}$ for $K \subseteq \mu^+$, 
are all distinct, so $B$ has $2^{\mu^+}$ ideals.) The argument is 
elementary
but a little tedious and we give it in some detail. Assume for contradiction 
that 
$b_i \leq \sum_{j \in J} b_j$. 

For arbitrary $n \in \omega$, we have the following situation. $d_{in}$ is 
non-zero
and for $j \in J$, $\{ d_{jm} : m \in \omega \}$ is a partition of unity; 
hence there
are elements $m(j) \in \omega$, for $j \in J$, such that $p = d_{in} \cdot 
\prod_{j \in J}
d_{jm(j)} $ is non-zero. Now $b_i \cdot d_{in} \leq z_{in} $ and thus $b_i
\cdot p \leq z_{in}$; similarly $b_j \cdot p \leq z_{jm(j)}$ holds for $j \in 
J$. It
follows from $b_i \leq \sum_{j \in J} b_j$ that $z_{in} \cdot p \leq b_i
\cdot p \leq \sum_{j \in J} z_{jm(j)}$. But $ \text{supp } p \subseteq X$ 
and 
 $z_{in}$, $z_{jm(j)}$ are in $Z$; hence $ z_{in} \leq \sum_{j \in J} 
z_{jm(j)}$. 
So $ z_{in} = z_{jm(j)}$, for some 
$j \in J$, since $Z \subseteq U$ is independent. Since $z_{in}$ was chosen 
in Step 7 from 
$Z_t$, where $t = f_i \restriction n$, and $(Z_t)_{t \in T}$ was a disjoint 
family,
it follows that $n = m(j)$ and $f_i \restriction n = f_j \restriction n$.

We have thus shown that for every $n \in \omega$, there is some $j \in J$ 
satisfying 
 $f_i \restriction n = f_j \restriction n$. But then $f_i \in \{f_j : j \in J \}$ 
and $i \in J$ (since the enumeration $\{ f_i : i < \mu^+ \}$ in Step 6 was
one-one), a contradiction.
\qed

\bigskip
\bigskip
\bigskip
\bigskip
\bigskip

\enddocument